\documentclass{article}

\usepackage{xparse}                                             

\usepackage{amsmath}
\usepackage{amsthm}
\usepackage{thmtools}

\numberwithin{equation}{section}
\numberwithin{figure}{section}
\numberwithin{table}{section}

\usepackage{amsfonts}
\usepackage{amssymb}
\usepackage{mathtools}
\makeatletter
\DeclareRobustCommand*\cal{\@fontswitch\relax\mathcal}
\makeatother

\usepackage{graphicx}
\usepackage{xcolor}

\usepackage[colorlinks]{hyperref}
\usepackage[capitalize]{cleveref}
\crefname{section}{Section}{Sections}
\Crefname{section}{Section}{Sections}
\crefname{construction}{Construction}{Constructions}
\Crefname{construction}{Construction}{Constructions}
\usepackage{url}
\newcommand{\doi}[1]{\textsc{doi}: \href{http://dx.doi.org/#1}{\nolinkurl{#1}}}

\usepackage{enumitem}
\setlist[enumerate, 1]{label = (\roman*), font = \upshape}
\setlist[itemize, 2]{label = {$\circ$}}
\usepackage{longtable}
\usepackage{subcaption}

\theoremstyle{plain}
\newtheorem{theorem}{Theorem}[section]
\newtheorem{corollary}[theorem]{Corollary}
\newtheorem{lemma}[theorem]{Lemma}

\newtheorem{proposition}[theorem]{Proposition}

\theoremstyle{definition}

\newtheorem{definition}[theorem]{Definition}
\newtheorem{example}[theorem]{Example}

\theoremstyle{remark}

\newcommand{\cA}{\mathcal{A}}
\newcommand{\cB}{\mathcal{B}}
\newcommand{\cC}{\mathcal{C}}

\newcommand{\cL}{\mathcal{L}}
\newcommand{\cF}{\mathcal{F}}
\newcommand{\cG}{\mathcal{G}}
\newcommand{\cH}{\mathcal{H}}

\DeclareMathOperator{\RR}{\mathbb{R}}

\DeclareMathOperator{\ZZ}{\mathbb{Z}}

\DeclareMathOperator{\supp}{\mathrm{supp}}

\newcommand{\from}{\colon}
\renewcommand{\mid}{:\ }
\newcommand{\gt}{>}
\newcommand{\lt}{<}

\title{Sizes of cross-intersecting families}
\date{}
\author{Ziqing Xiang\thanks{Southern University of Science and Technology. Email: \texttt{xiangzq@sustech.edu.cn}.}
\and Huajun Zhang\thanks{Shaoxing University. Email: \texttt{huajunzhang@usx.edu.cn}.}}

\begin{document}

\maketitle

\begin{abstract}
The size of cross-intersecting families $\cF_1, \dots, \cF_d$ is usually measured by the sum or the product of the sizes $|\cF_i|$. These two natural size functions have different extremal structures, and different techniques have been developed for them.

We introduce the cross-degree of a size function and show that the potential extremal structures of maximum cross-intersecting families are the same for all size functions of finite cross-degree, including the sum and the product. Our main tool is a strengthened form of the generating family method of Ahlswede--Khachatrian, based on standard generating families and on local modifications of cross-intersecting families, which we call mutations. It treats two completely different types of extremal structures in a unified way and, unlike the original method, needs shifting in only one step.
\end{abstract}

\section{Introduction}

Families $\cF_1, \dots, \cF_d$ are {\em cross-intersecting} if $A \cap B \neq \emptyset$ for all $A \in \cF_i$ and $B \in \cF_j$ with $i \neq j$. They generalize intersecting families, which were first studied by Erd\H{o}s--Ko--Rado \cite{EKR61} in 1961.

In 1967, Hilton--Milner \cite{HM67} initiated the study of cross-intersecting families whose sizes are measured by the sum $|\cF_1| + \dots + |\cF_d|$. They showed that, for $d = 2$ and $n \geq 2 k$, the maximum of the sum for nonempty families is attained by $\{S\}$ and the family of all $k$-sets meeting $S$, for some $k$-set $S$. The sum has also been studied for more than two families, possibly empty \cite{Bor09, Hil77} or nonempty \cite{HP25, SFQ22, ZF24}, and for families of different uniformities \cite{FT92, HP25, SFQ22, ZF24}. Further results concern weighted sums \cite{HP25, SFQ22}, nontrivial families, whose members have no common element \cite{Fra24}, and cross-$t$-intersecting families \cite{WZ11, WZ13}. Analogous results have been obtained for permutations \cite{Bor10}, signed sets \cite{BL10}, labeled sets \cite{Bor08} and vector spaces \cite{WZ11, WZ13}, and, when $d$ is large, for subfamilies of an arbitrary family \cite{Bor14a}.

In 1986, Pyber \cite{Pyb86} initiated the study of cross-intersecting families whose sizes are measured by the product $|\cF_1| \cdots |\cF_d|$. He showed that, for two families and $n \geq 2 k$, the product is maximized by two copies of a star. The product has also been studied for more than two families \cite{Bor14b, Bor17, FT11}, for families of different uniformities \cite{Bor14b, MT89}, for nontrivial families \cite{FK17, FW23, Hua26}, and for cross-$t$-intersecting families \cite{Bor14b, FLST14, Tok13}. Analogous results have been obtained for permutations \cite{Bor17} and vector spaces \cite{ST14}, and, when $d$ is large, for subfamilies of an arbitrary family \cite{Bor14a}.

Both size functions are natural, with no compelling reason to prefer either. However, for fixed $d$ and large $n$, their extremal structures are completely different: for the sum, all families but one are nearly empty, while for the product, all families are close to stars. As a result, different techniques have been developed for them. This dichotomy calls for a unified explanation, and the key difficulty is to handle at least these two types of extremal structures at the same time.

In this paper, we provide such an explanation by developing tools that overcome this difficulty. A fundamental notion is the cross-degree of a size function. For example, the sum $x_1 + \dots + x_d$ has cross-degree at most $0$, and the product $x_1 \cdots x_d$ has cross-degree at most $1$, in the sense of the following definition.

\begin{definition} \label{def:cross-degree}
Let $c$ be a nonnegative real number, and let $f$ be a function of positive real variables $x_1, \dots, x_d$. We say that {\em $f$ has cross-degree at most $c$} if
\begin{enumerate}
  \item \label{itm:twice-differentiable} $f$ is twice continuously differentiable,
  \item \label{itm:derivative-signs} $\partial_{x_i} f \gt 0$ and $\partial_{x_i x_j} f \geq 0$ for $1 \leq i, j \leq d$,
  \item \label{itm:cross-condition} $c \ \partial_{x_i} f \geq x_j \ \partial_{x_i x_j} f$ for $1 \leq i, j \leq d$ with $i \neq j$.
\end{enumerate}
\end{definition}

The main result of this paper characterizes the maximum cross-intersecting families with respect to an arbitrary size function of finite cross-degree.

\begin{theorem} \label{thm:main}
Let $f \from \RR_{\gt 0}^d \to \RR$ be a function of cross-degree at most $c$ with $d \geq 2$, where $c = 0$ or $c \geq 1$. Let $\cF_1, \dots, \cF_d \subseteq \binom{[n]}{k}$ be nonempty cross-intersecting families that maximize $f(|\cF_1|, \dots, |\cF_d|)$ among all nonempty cross-intersecting families. Suppose that $n \geq (k - 1)(c + 2) + 3$. Then $(\cF_1, \dots, \cF_d)$ is a permutation of $(\cA^S, \dots, \cA^S, \cB^S)$ for some nonempty $S \subseteq [n]$ with $|S| \leq k$, where
\begin{align*}
  \cA^S & := \left\{ A \in \binom{[n]}{k} \mid S \subseteq A \right\}, \\
  \cB^S & := \left\{ B \in \binom{[n]}{k} \mid B \cap S \neq \emptyset \right\}.
\end{align*}
Moreover, the following statements hold.
\begin{enumerate}
\item \label{itm:main-at-most-c} If $n \geq 2 (c + 1) k$, then either $|S| = 1$, or $|S| = k$, or $2 \leq |S| \leq \lceil c \rceil$.
\item \label{itm:main-one-or-k} If $n \geq (k - 1) 2^{c + 1} + 3$, then either $|S| = 1$, or $|S| = k$.
\end{enumerate}
\end{theorem}

In particular, \cref{thm:main} covers both types of extremal structures: copies of a star for $|S| = 1$, and, for $|S| = k$, the structure of Hilton--Milner, in which all families but one equal $\{S\}$. Note that the assumption on $n$ in \cref{thm:main} cannot be weakened to $n \geq 2 k + 1$ for general $c$, as shown in \cref{ex:small-n}.

We give applications of \cref{thm:main} in \cref{sec:examples}, explain the main ideas of its proof in \cref{sec:main-ideas}, and prove it in \cref{sec:proof}.

\subsection{Applications and examples} \label{sec:examples}

In this subsection, we apply \cref{thm:main} to various size functions. We recover the classical results for the sum and the product, find extremal structures that lie between the two types, and show that the assumption on $n$ in \cref{thm:main} cannot be weakened to $n \geq 2 k + 1$.

\begin{example} \label{ex:sum}
Consider the classical sum size function $f := x_1 + \dots + x_d$, which has cross-degree at most $0$. Hence, by \cref{thm:main}\ref{itm:main-one-or-k}, for $n \geq 2 k + 1$, the only potential extremal structures have $|S| = 1$ or $|S| = k$. A direct comparison shows that the maximum of $|\cF_1| + \dots + |\cF_d|$ is
\[
  \max \left\{ d \binom{n - 1}{k - 1}, \binom{n}{k} - \binom{n - k}{k} + d - 1 \right\}.
\]
\end{example}

For $d = 2$, \cref{ex:sum} recovers the result of Hilton--Milner \cite{HM67} for $n \geq 2 k + 1$, while \cite{HM67} also covers $n = 2 k$. For $d \geq 3$, it recovers the result of Shi--Frankl--Qian \cite{SFQ22} for $n \geq 2 k + 1$, while \cite{SFQ22} also covers $n = 2 k$.

A variant of \cref{thm:main} without the nonemptiness assumption is easily obtained, provided that the size function extends to $\RR_{\geq 0}^d$ and its restrictions to fewer variables still have cross-degree at most $c$, as is the case for the sum size function. Indeed, it suffices to apply \cref{thm:main} to these restrictions, and the case of a single nonempty family is trivial. Using this variant, we also recover Hilton's theorem \cite{Hil77} for the same range of $n$.

\begin{example} \label{ex:product}
Consider the classical product size function $f := x_1 \cdots x_d$, which has cross-degree at most $1$. Hence, by \cref{thm:main}, for $n \geq 3 k$, the only potential extremal structures consist of $d - 1$ copies of $\cA^S$ and one copy of $\cB^S$. A direct comparison shows that only $|S| = 1$ gives the maximum, so the extremal structure consists of $d$ copies of a star.
\end{example}

For $d = 2$, \cref{ex:product} recovers the result of Pyber \cite{Pyb86} for $n \geq 3 k$, while \cite{Pyb86} covers all $n \geq 2 k$. For $d \geq 3$, the conclusion can also be deduced from the case $d = 2$ by considering the families in pairs.

\begin{example} \label{ex:interpolation}
Consider $f := (x_1 + \gamma)(x_2 + \gamma) = x_1 x_2 + \gamma (x_1 + x_2) + \gamma^2$, where $\gamma \geq 0$ is a real number. This function interpolates between the product and the sum, and it has cross-degree at most $1$. Let $k \geq 3$. By \cref{thm:main}\ref{itm:main-at-most-c}, if $n \geq 4 k$, then the only potential extremal structures have $|S| = 1$ or $|S| = k$. A direct comparison shows that the maximum of $f(|\cF_1|, |\cF_2|)$ is attained by two stars if $\gamma \lt \gamma^*$, and by $\{S\}$ and $\cB^S$ with $|S| = k$ if $\gamma \gt \gamma^*$. Here $\gamma^* \gt 0$ is the unique value of $\gamma$ for which both give the same value of $f$.
\end{example}

\cref{ex:interpolation} shows that the extremal structures for the product and for the sum, which were obtained by different methods, are two phases of a single one-parameter family of size functions.

\begin{example} \label{ex:cross-degree-5}
Consider $f := x_1 x_2^5$, which has cross-degree at most $5$. Let $k = 4$. Then \cref{thm:main} applies for $n \geq 24$, and a direct computation, using \cref{thm:main}\ref{itm:main-one-or-k} for $n \geq 195$, shows that the extremal structure is $(\cF_1, \cF_2) = (\cA^S, \cB^S)$ with
\begin{itemize}
\item $|S| = 2$ if $24 \leq n \leq 89$,
\item $|S| = 1$, namely a pair of stars, if $n \geq 90$.
\end{itemize}
\end{example}

\begin{example} \label{ex:cross-degree-10}
Consider $f := x_1 x_2^{10}$, which has cross-degree at most $10$. Let $k = 4$. Then \cref{thm:main} applies for $n \geq 39$, and a direct computation, using \cref{thm:main}\ref{itm:main-one-or-k} for $n \geq 6147$, shows that the extremal structure is $(\cF_1, \cF_2) = (\cA^S, \cB^S)$ with
\begin{itemize}
\item $|S| = 3$ if $39 \leq n \leq 101$,
\item $|S| = 2$ if $102 \leq n \leq 3057$,
\item $|S| = 1$, namely a pair of stars, if $n \geq 3058$.
\end{itemize}
\end{example}

\cref{ex:cross-degree-5,ex:cross-degree-10} exhibit extremal structures with $2 \leq |S| \leq k - 1$, which occur neither for the sum nor for the product, and in \cref{ex:cross-degree-10} all values $2 \leq |S| \leq k - 1$ occur.

\begin{example} \label{ex:polynomial}
Consider
\[
  f := \sum_{\kappa \in K} \lambda_{\kappa} x_1^{\kappa_1} \cdots x_d^{\kappa_d}
\]
with $\lambda_{\kappa} \gt 0$, where $K \subseteq (\{0\} \cup \{ t \in \RR \mid t \geq 1 \})^d$ is a finite set containing, for each $i$, an element $\kappa$ with $\kappa_i \neq 0$. A direct calculation shows that $f$ has cross-degree at most $c$, where $c := \max \{ \kappa_i \mid \kappa \in K, 1 \leq i \leq d \}$. Since $c \geq 1$, each of the following assumptions on $n$ implies $n \geq (k - 1)(c + 2) + 3$. If $n \geq 2 (c + 1) k$, then \cref{thm:main}\ref{itm:main-at-most-c} shows that, up to a permutation of $[n]$, there are at most $d \lceil c \rceil + 1$ potential extremal structures. If $n \geq (k - 1) 2^{c + 1} + 3$, then \cref{thm:main}\ref{itm:main-one-or-k} shows that there are at most $d + 1$ potential extremal structures.
\end{example}

\begin{example} \label{ex:small-n}
Consider $f := x_1 x_2^2$, which has cross-degree at most $2$. Let $k = 4$ and $n = 9$. Then the families
\begin{align*}
  \cA & := \left\{ A \in \binom{[9]}{4} \mid 1 \in A, A \cap \{2, 3\} \neq \emptyset \right\}, \\
  \cB & := \left\{ B \in \binom{[9]}{4} \mid 1 \in B \text{ or } \{2, 3\} \subseteq B \right\}
\end{align*}
are cross-intersecting with $|\cA| = 36$ and $|\cB| = 71$, and $f(36, 71) = 181476$ exceeds the values of $f$ at the sizes of all pairs $(\cA^S, \cB^S)$ and $(\cB^S, \cA^S)$, the largest of which is $f(56, 56) = 175616$.
\end{example}

\subsection{Main ideas} \label{sec:main-ideas}

To handle both types of extremal structures at the same time, we compare a cross-intersecting pair with two local modifications of it in opposite directions. Our approach builds on the generating family method, which was introduced by Ahlswede--Khachatrian \cite{AK97} in 1997 to prove the complete intersection theorem. In this method, a left-shifted family is described by its generating sets, it is modified locally at the generating sets that contain the largest element, and two such modifications are compared.

Our standard generating families (see \cref{def:standard-generating-family}), mutations (see \cref{def:mutation}) and hinges (see \cref{def:hinge}) are strengthened forms of these three ingredients. They apply to families that are not necessarily left-shifted, and they show more clearly why the method works. As a result, our proof uses shifting in only one step, namely to show that hinges exist.

Every $k$-uniform family has a canonical generating family, the standard generating family, which has the smallest support among all generating families, see \cref{sec:generating}. For a maximal cross-intersecting pair, the two standard generating families have the same support. A mutation of the pair is determined by a member of either standard generating family. It shrinks one family and enlarges the other, and the result is again cross-intersecting, see \cref{sec:mutations}. A hinge is a pair of members, one from each standard generating family, that cover the common support and meet in exactly one element. It gives two mutations in opposite directions, see \cref{sec:hinges}.

Mutations at hinges explain why the potential extremal structures in \cref{thm:main} are the same for all size functions. For $d = 2$, these structures are the pairs $(\cA^S, \cB^S)$ and $(\cB^S, \cA^S)$, which we call {\em distinguished}. Unless a pair with a hinge is distinguished, the two mutations at the hinge consist of nonempty families, see \cref{prop:hinge-cases}.

For every size function of finite cross-degree and $n$ large enough, one of the two mutations increases the value of the size function, unless both families have fewer members than a star, see \cref{thm:hinge-comparison}. Apart from monotonicity, the size function enters this comparison only through \cref{lem:three-points}, which compares its values at three points. Since two stars with the same center are cross-intersecting, it follows that a maximum pair with a hinge is distinguished.

The distinguished pairs $(\cA^S, \cB^S)$ are themselves related by mutations. Up to a permutation of $[n]$, they form a chain that starts with the pair of stars, the extremal structure for the product, and ends with an extremal structure for the sum of two nonempty families, see \cref{sec:chain}. Comparing three consecutive pairs in this chain, we obtain the bounds on $|S|$ in \cref{thm:main}\ref{itm:main-at-most-c} and \ref{itm:main-one-or-k}, see \cref{prop:consecutive}.

Hinges need not exist in general, see \cref{sec:hinges}. To obtain them, we replace the families by initial segments of the lexicographic order of the same sizes, which are cross-intersecting by Hilton's lemma. Initial segments are left-shifted, and maximal cross-intersecting pairs of nonempty left-shifted families have hinges, see \cref{prop:hinge-exists}. We then recover the original families from their sizes, using the known fact that distinguished pairs are determined by their sizes, see \cref{sec:determined-by-sizes}.

Functions of bounded cross-degree are studied in \cref{sec:size-functions}, standard generating families in \cref{sec:generating}, mutations and hinges in \cref{sec:mutations}, maximum pairs with a hinge in \cref{sec:maximum}, distinguished pairs in \cref{sec:distinguished}, and hinges of left-shifted families in \cref{sec:left-shifted}. Finally, \cref{thm:main} is proved in \cref{sec:proof}.

\subsection*{Notation}

Throughout the paper, $d$, $n$ and $k$ denote positive integers. Unless stated otherwise, we consider only cross-intersecting families in $\binom{[n]}{k}$. We assume that $d \geq 2$ and $n \geq 2 k$, since a single family is trivially cross-intersecting, and any two $k$-subsets of $[n]$ intersect if $n \lt 2 k$. For integers $a$ and $b$, let $[a, b] := \{ i \in \ZZ \mid a \leq i \leq b \}$ and $[b] := [1, b]$. For a set $U$ and an integer $r$, we write $\binom{U}{r}$, $\binom{U}{\leq r}$ and $\binom{U}{\geq r}$ for the families of subsets of $U$ with exactly $r$, at most $r$ and at least $r$ elements, respectively. A binomial coefficient with a negative lower index is $0$.

\section{Size functions of bounded cross-degree} \label{sec:size-functions}

In this section, we collect the properties of functions of bounded cross-degree that are needed later. In the proof of \cref{thm:main}, apart from monotonicity, the conditions in \cref{def:cross-degree} are used only through two lemmas. \cref{lem:two-variables} reduces functions of $d$ variables to functions of two variables, and \cref{lem:three-points} compares the values of a function of two variables at three points.

For functions of two variables, we allow different constants for the two variables in condition \ref{itm:cross-condition} of \cref{def:cross-degree}. This keeps the bounds in \cref{thm:main} independent of $d$, since the reduction from $d$ variables in \cref{lem:two-variables} multiplies only one of the two constants by $d - 1$.

\begin{definition} \label{def:cross-degrees}
Let $c_1$ and $c_2$ be nonnegative real numbers. We say that a function $f(x, y)$ of two positive real variables has {\em cross-degrees at most $(c_1, c_2)$} if $f$ is twice continuously differentiable, satisfies \cref{def:cross-degree}\ref{itm:derivative-signs}, that is,
\[
  \partial_x f \gt 0, \quad \partial_y f \gt 0, \quad \partial_{xx} f, \partial_{xy} f, \partial_{yy} f \geq 0,
\]
and satisfies
\[
  c_1 \ \partial_x f \geq y \ \partial_{xy} f, \quad c_2 \ \partial_y f \geq x \ \partial_{xy} f.
\]
\end{definition}

Thus a function of two variables has cross-degree at most $c$ if and only if it has cross-degrees at most $(c, c)$. The following lemma, which is used in \cref{sec:proof}, reduces functions of $d$ variables to functions of two variables.

\begin{lemma} \label{lem:two-variables}
Let $f$ be a function of $d$ variables with cross-degree at most $c$. For $p \in [d]$, let
\[
  g_p(x, y) := f(x, \dots, x, y, x, \dots, x),
\]
where $y$ is in the $p$-th position. Then $g_p$ has cross-degrees at most $(c, (d - 1) c)$.
\end{lemma}

\begin{proof}
Let $I := [d] \setminus \{p\}$. The partial derivatives of $g_p$ of order at most two are sums of partial derivatives of $f$ at $(x, \dots, x, y, x, \dots, x)$, for example,
\begin{align*}
  \partial_x g_p & = \sum_{j \in I} \partial_{x_j} f, & \partial_y g_p & = \partial_{x_p} f, \\
  \partial_{x x} g_p & = \sum_{j, l \in I} \partial_{x_j x_l} f, & \partial_{x y} g_p & = \sum_{j \in I} \partial_{x_j x_p} f.
\end{align*}
Hence $g_p$ is twice continuously differentiable and satisfies \cref{def:cross-degree}\ref{itm:derivative-signs}. Moreover, \cref{def:cross-degree}\ref{itm:cross-condition} for $f$ gives $c \ \partial_x g_p \geq y \ \partial_{x y} g_p$ and $(d - 1) c \ \partial_y g_p \geq x \ \partial_{x y} g_p$, since $|I| = d - 1$.
\end{proof}

The following lemma bounds both coordinates of a point at which the value of a function does not increase under two moves in opposite directions. It is used to compare the value of a size function at a pair with its values at the two mutations at a hinge in \cref{thm:hinge-comparison}, and at the two neighboring distinguished pairs in \cref{prop:consecutive}.

\begin{lemma} \label{lem:three-points}
Let $f(x, y)$ be a function of cross-degrees at most $(c_1, c_2)$, where $c_i = 0$ or $c_i \geq 1$ for $i = 1, 2$. Let $a, b, a_1, b_1, a_2, b_2$ be positive real numbers with $a \gt a_2$, $b \gt b_1$ and $\frac{b_2}{b_1} \gt \frac{a_2}{a_1}$. If
\[
  f(a, b) \geq f(a + a_1, b - b_1), \quad \text{and} \quad f(a, b) \geq f(a - a_2, b + b_2),
\]
then
\begin{align}
  \left(\frac{b}{b - b_1}\right)^{c_1} & \geq 1 + \frac{a_1 b_2 - a_2 b_1}{a_2 (b_1 + b_2)}, \label{eq:ratio-b} \\
  \left(\frac{a}{a - a_2}\right)^{c_2} & \geq 1 + \frac{a_1 b_2 - a_2 b_1}{b_1 (a_1 + a_2)}. \label{eq:ratio-a}
\end{align}
In particular, $c_1 \geq 1$ and $c_2 \geq 1$, and
\begin{equation} \label{eq:bound-a}
  a \leq c_2 \left(\frac{b_2}{b_1} - \frac{a_2}{a_1}\right)^{-1} \left(1 + \frac{b_2}{b_1}\right) a_2,
\end{equation}
and
\begin{equation} \label{eq:bound-b}
  b \leq c_1 \left(\frac{b_2}{b_1} - \frac{a_2}{a_1}\right)^{-1} \left(1 + \frac{a_2}{a_1}\right) b_2.
\end{equation}
\end{lemma}

\begin{proof}
We compare the increments of $f$ along horizontal and vertical segments near $(a, b)$. Let $\rho := \frac{b}{b - b_1}$. Since $\partial_x f \gt 0$, the condition $c_1 \ \partial_x f \geq y \ \partial_{x y} f$ means that $\partial_y \log (\partial_x f) \leq \frac{c_1}{y}$. Integrating this inequality from $y = b - b_1$ to $y = b$, we obtain
\begin{equation} \label{eq:partial-ratio}
  (\partial_{x} f)(x, b) \leq \rho^{c_1} (\partial_{x} f)(x, b - b_1)
\end{equation}
for all $x \gt 0$. Let
\begin{align*}
  H_1 & := f(a + a_1, b - b_1) - f(a, b - b_1), & H_2 & := f(a, b - b_1) - f(a - a_2, b - b_1), \\
  V_1 & := f(a - a_2, b) - f(a - a_2, b - b_1), & V_2 & := f(a - a_2, b + b_2) - f(a - a_2, b),
\end{align*}
and
\[
  M := f(a, b) - f(a - a_2, b) - f(a, b - b_1) + f(a - a_2, b - b_1).
\]
Then $f(a, b) - f(a, b - b_1) = V_1 + M$ and $f(a, b) - f(a - a_2, b) = H_2 + M$. Hence the two assumptions say that
\[
  H_1 \leq V_1 + M, \quad \text{and} \quad V_2 \leq H_2 + M.
\]
Since $\partial_{x x} f \geq 0$ and $\partial_{y y} f \geq 0$, we have $a_2 H_1 \geq a_1 H_2$ and $b_1 V_2 \geq b_2 V_1$. Moreover, integrating \cref{eq:partial-ratio} from $x = a - a_2$ to $x = a$, we obtain
\[
  H_2 + M = \int_{a - a_2}^{a} (\partial_{x} f)(x, b) \, dx \leq \rho^{c_1} \int_{a - a_2}^{a} (\partial_{x} f)(x, b - b_1) \, dx = \rho^{c_1} H_2.
\]
Therefore,
\[
  a_1 b_2 H_2 \leq a_2 b_2 H_1 \leq a_2 b_2 (V_1 + M) \leq a_2 b_1 V_2 + a_2 b_2 M \leq a_2 b_1 (H_2 + M) + a_2 b_2 M,
\]
that is,
\[
  (a_1 b_2 - a_2 b_1) H_2 \leq a_2 (b_1 + b_2) M \leq a_2 (b_1 + b_2) (\rho^{c_1} - 1) H_2.
\]
Since $\partial_x f \gt 0$, we have $H_2 \gt 0$, and \cref{eq:ratio-b} follows.

The function $\overline{f}(x, y) := f(y, x)$ has cross-degrees at most $(c_2, c_1)$. Moreover, the assumptions for $f$ and $(a, b, a_1, b_1, a_2, b_2)$ are exactly the assumptions for $\overline{f}$ and $(b, a, b_2, a_2, b_1, a_1)$. Hence \cref{eq:ratio-a} follows from \cref{eq:ratio-b} applied to $\overline{f}$.

Since $\frac{b_2}{b_1} \gt \frac{a_2}{a_1}$, the right-hand sides of \cref{eq:ratio-b} and \cref{eq:ratio-a} are larger than $1$. Hence $c_1$ and $c_2$ are nonzero, that is, $c_1 \geq 1$ and $c_2 \geq 1$. Let
\[
  t := \frac{b_1}{b} \quad \text{and} \quad q := \frac{a_1 b_2 - a_2 b_1}{a_2 (b_1 + b_2)}.
\]
By Bernoulli's inequality and \cref{eq:ratio-b},
\[
  1 - c_1 t \leq (1 - t)^{c_1} = \rho^{- c_1} \leq \frac{1}{1 + q},
\]
hence $c_1 t \geq \frac{q}{1 + q}$, that is, $b \leq c_1 \frac{1 + q}{q} b_1$. A direct calculation shows that
\[
  \frac{1 + q}{q} b_1 = \frac{(a_1 + a_2) b_1 b_2}{a_1 b_2 - a_2 b_1} = \left(\frac{b_2}{b_1} - \frac{a_2}{a_1}\right)^{-1} \left(1 + \frac{a_2}{a_1}\right) b_2,
\]
which proves \cref{eq:bound-b}. In the same way, \cref{eq:bound-a} follows from \cref{eq:ratio-a}.
\end{proof}

\section{Standard generating families} \label{sec:generating}

The generating family method of Ahlswede--Khachatrian \cite{AK97} plays a fundamental role in extremal set theory, but in their formulation and in the book of Gerbner--Patk\'os \cite{GP18}, it applies only to left-shifted families. In this section and in \cref{sec:mutations}, we generalize it to families that are not necessarily left-shifted. Every $k$-uniform family has a canonical generating family, which we call the standard generating family. It has the smallest support among all generating families. For a maximal cross-intersecting pair, the two standard generating families have the same support, see \cref{thm:generating}. The mutations in \cref{sec:mutations} are defined in terms of this common support.

Throughout this section, we fix a ground set $U$, and all families are subfamilies of $2^{U}$. We state explicitly when a family is $k$-uniform.

\subsection{Operations on families} \label{sec:operations}

In this subsection, we introduce the support, the up-closure, the trace and the restriction.

The {\em support} of a family $\cF$ is its union, that is, $\supp \cF := \bigcup_{F \in \cF} F$. For $X \subseteq U$, let $X^\uparrow := \{Y \subseteq U \mid X \subseteq Y\}$. The {\em up-closure} of a family $\cF$ is
\[
  \cF^\uparrow := \bigcup_{F \in \cF} F^\uparrow = \{ G \subseteq U \mid F \subseteq G \text{ for some } F \in \cF \}.
\]

For a subset $S \subseteq U$, we call $Y \cap S$ the {\em trace} of $Y \subseteq U$ on $S$. For $X \subseteq S$, let
\[
  \langle X \rangle_S := \left\{ F \in \binom{U}{k} \mid F \cap S = X \right\}
\]
be the family of $k$-sets whose trace on $S$ is $X$. Clearly, $|\langle X \rangle_S| = \binom{|U| - |S|}{k - |X|}$.

For a subset $X \subseteq U$, the {\em restriction} of $\cF$ to $X$ is
\[
  \cF|_X := \{ X \cap F \mid F \in \cF \}.
\]

\subsection{The standard generating family}

Consider a $k$-uniform family $\cF$. A {\em generating family} of $\cF$ is a family $\cG \subseteq \binom{U}{\leq k}$ such that $\cF = \cG^\uparrow \cap \binom{U}{k}$, or equivalently, $\cF^\uparrow = \cG^\uparrow \cap \binom{U}{\geq k}$. Among the generating families of $\cF$, there is a canonical one.

\begin{definition} \label{def:standard-generating-family}
The {\em standard generating family of $\cF$}, denoted by $\cF_{\rm std}$, is the family of inclusion-minimal members of
\begin{equation} \label{eq:std-family}
\left\{ G \subseteq U \mid G^\uparrow \cap \binom{U}{\geq k} \subseteq \cF^\uparrow \right\}.
\end{equation}
\end{definition}

Note that the family in \cref{eq:std-family} is closed under taking supersets, hence it coincides with $\cF_{\rm std}^\uparrow$. For $G \subseteq U$ with $|G| \leq k \leq |U|$, we have $G \in \cF_{\rm std}^\uparrow$ if and only if every $k$-subset of $U$ containing $G$ belongs to $\cF$. Every member of $\cF_{\rm std}$ has at most $k$ elements: if $G \in \cF_{\rm std}^\uparrow$ with $|G| \gt k$, then $G$ contains some $F \in \cF$, and $F \in \cF_{\rm std}^\uparrow$, so $G$ is not minimal. It is straightforward to verify that $\cF_{\rm std}^\uparrow \cap \binom{U}{\leq k}$ is a generating family of $\cF$. Moreover, every generating family $\cG$ of $\cF$ is a subfamily of $\cF_{\rm std}^\uparrow \cap \binom{U}{\leq k}$.

The standard generating family $\cF_{\rm std}$ is itself a generating family, since it is the family of inclusion-minimal members of the generating family $\cF_{\rm std}^\uparrow \cap \binom{U}{\leq k}$. In particular, every $F \in \cF$ contains some $G \in \cF_{\rm std}$, and $G \subseteq F \cap \supp \cF_{\rm std}$. Note also that $\cF_{\rm std}$ is an antichain. The standard generating family has the smallest support among all generating families.

\begin{proposition} \label{prop:std-support}
Let $\cF \subseteq \binom{U}{k}$ be a family, and $\cG$ a generating family of $\cF$. Then $\supp \cF_{\rm std} \subseteq \supp \cG$.
\end{proposition}

\begin{proof}
Let $H \in \cF_{\rm std}$, and assume to the contrary that there exists $h \in H \setminus \supp \cG$. We claim that $H \setminus \{h\} \in \cF_{\rm std}^\uparrow$, which contradicts the minimality of $H$. Let $F \supseteq H \setminus \{h\}$ be a $k$-set. If $h \in F$, then $F \supseteq H$, hence $F \in \cF$. If $h \notin F$, take $w \in F \setminus H$, which exists since $|F| = k \geq |H| \gt |H \setminus \{h\}|$, and let $F' := (F \setminus \{w\}) \cup \{h\}$. Then $F' \supseteq H$, hence $F' \in \cF$ and $F'$ contains some $G \in \cG$. Since $h \notin \supp \cG$, we have $G \subseteq F' \setminus \{h\} \subseteq F$, hence $F \in \cF$. This proves the claim.
\end{proof}

\subsection{Cross-intersecting families}

Two families $\cF$ and $\cF'$ are {\em cross-intersecting} if $F \cap F' \neq \emptyset$ for all $F \in \cF$ and $F' \in \cF'$. For families $\cF_1, \dots, \cF_d$, we say that $(\cF_1, \dots, \cF_d)$ is {\em cross-intersecting} if $\cF_i$ and $\cF_j$ are cross-intersecting for all distinct $i$ and $j$. A cross-intersecting tuple $(\cF_1, \dots, \cF_d)$ of $k$-uniform families is {\em maximal} if, for every cross-intersecting $(\cF_1', \dots, \cF_d')$ with $\cF_i \subseteq \cF_i' \subseteq \binom{U}{k}$ for all $i$, we have $\cF_i' = \cF_i$ for all $i$.

For a family $\cF$ of $k$-subsets of $U$, let
\[
  \cF^* := \left\{ X \in \binom{U}{k} \mid X \cap F \neq \emptyset \text{ for all } F \in \cF \right\}
\]
be the largest $k$-uniform family that is cross-intersecting with $\cF$. A pair $(\cA, \cB)$ of $k$-uniform families is maximal cross-intersecting if and only if $\cA = \cB^*$ and $\cB = \cA^*$. More generally, cross-intersecting $k$-uniform families $(\cF_1, \dots, \cF_d)$ are maximal if and only if $\cF_i = (\bigcup_{j \neq i} \cF_j)^*$ for all $i$.

The following theorem shows that cross-intersection can be checked on generating families, and that the standard generating families of a maximal cross-intersecting pair have the same support.

\begin{theorem} \label{thm:generating}
Suppose that $|U| \geq 2 k$. Let $\cA, \cB$ be two $k$-uniform families.
\begin{enumerate}
\item \label{itm:cross-via-generating} Let $\cG, \cH$ be generating families of $\cA, \cB$, respectively. Then $(\cA, \cB)$ is cross-intersecting if and only if $(\cG, \cH)$ is cross-intersecting.
\item \label{itm:common-support} If $(\cA, \cB)$ is maximal cross-intersecting, then $\supp \cA_{\rm std} = \supp \cB_{\rm std}$.
\end{enumerate}
\end{theorem}

\begin{proof}
\ref{itm:cross-via-generating}
Since $\cA \subseteq \cG^\uparrow$ and $\cB \subseteq \cH^\uparrow$, if $(\cG, \cH)$ is cross-intersecting, then so is $(\cA, \cB)$.

Conversely, assume that $(\cA, \cB)$ is cross-intersecting while $(\cG, \cH)$ is not. Take $G \in \cG$ and $H \in \cH$ such that $G \cap H = \emptyset$. Since $|U| \geq 2k$, there exist a $(k - |G|)$-set $X \subseteq U \setminus (G \cup H)$ and a $(k - |H|)$-set $Y \subseteq U \setminus (G \cup H \cup X)$. Then $G \cup X \in \cA$ and $H \cup Y \in \cB$ are disjoint, which is a contradiction.

\ref{itm:common-support} Let $Z := (\supp \cA_{\rm std}) \cap (\supp \cB_{\rm std})$. For $G \in \cA_{\rm std}$ and $H \in \cB_{\rm std}$, we have $G \cap H \neq \emptyset$ by \ref{itm:cross-via-generating}, and clearly $G \cap H \subseteq Z$. Therefore, $\cA_{\rm std}|_Z$ and $\cB_{\rm std}|_Z$ are cross-intersecting. Let $\cA' := (\cA_{\rm std}|_Z)^\uparrow \cap \binom{U}{k}$ and $\cB' := (\cB_{\rm std}|_Z)^\uparrow \cap \binom{U}{k}$ be the $k$-uniform families generated by $\cA_{\rm std}|_Z$ and $\cB_{\rm std}|_Z$. By \ref{itm:cross-via-generating}, $(\cA', \cB')$ is cross-intersecting. Since $\cA \subseteq \cA'$, $\cB \subseteq \cB'$ and $(\cA, \cB)$ is maximal cross-intersecting, we have $(\cA, \cB) = (\cA', \cB')$. Hence $\cA_{\rm std}|_Z$ is a generating family of $\cA$, and $\supp \cA_{\rm std} \subseteq Z$ by \cref{prop:std-support}. Similarly, $\supp \cB_{\rm std} \subseteq Z$, which forces $\supp \cA_{\rm std} = Z = \supp \cB_{\rm std}$.
\end{proof}

\section{Mutations of cross-intersecting families} \label{sec:mutations}

In this section, we introduce {\em mutations}, which are local modifications of a maximal cross-intersecting pair. A mutation is determined by a member $X$ of one of the two standard generating families. It removes from one family all members whose trace on the common support is $X$, and adds to the other family all $k$-sets whose trace is the complement of $X$ in the common support. The size $|X|$ is called the rank of the mutation.

We first define mutations and discuss their basic properties. Then, we introduce hinges, which give pairs of mutations in opposite directions. Finally, we show that the two mutations at a hinge consist of nonempty families unless the pair is distinguished, that is, of the form $(\cA^S, \cB^S)$ or $(\cB^S, \cA^S)$.

Throughout this section, $U = [n]$ with $n \geq 2 k$, and we fix a maximal cross-intersecting pair $(\cA, \cB)$ of $k$-uniform families. By \cref{thm:generating}\ref{itm:common-support}, $\supp \cA_{\rm std} = \supp \cB_{\rm std}$, and we denote this common support by $Z$. For $X \subseteq Z$, we write $\langle X \rangle := \langle X \rangle_Z$ for the family of $k$-sets whose trace on $Z$ is $X$, see \cref{sec:operations}.

\begin{definition} \label{def:mutation}
For $X \in \cA_{\rm std}$, the {\em first mutation of $(\cA, \cB)$ at $X$} is
\[
  \mu^1_X(\cA, \cB) := \left(\cA \setminus \langle X \rangle, \ \cB \cup \langle Z \setminus X \rangle \right).
\]
For $Y \in \cB_{\rm std}$, the {\em second mutation of $(\cA, \cB)$ at $Y$} is
\[
  \mu^2_Y(\cA, \cB) := \left(\cA \cup \langle Z \setminus Y \rangle, \ \cB \setminus \langle Y \rangle \right).
\]
The {\em rank} of $\mu^1_X$ is $|X|$, and the {\em rank} of $\mu^2_Y$ is $|Y|$.
\end{definition}

Mutations of a maximal cross-intersecting pair are cross-intersecting, and their sizes can be computed explicitly.

\begin{proposition} \label{prop:mutation-properties}
Let $X \in \cA_{\rm std}$. Then the following statements hold.
\begin{enumerate}
\item \label{itm:mutation-traces} $\langle X \rangle \subseteq \cA$ and $\langle Z \setminus X \rangle \cap \cB = \emptyset$.
\item \label{itm:mutation-cross-intersecting} The first mutation $\mu^1_X(\cA, \cB)$ is cross-intersecting.
\item \label{itm:mutation-sizes} $|\langle X \rangle| = \binom{n - |Z|}{k - |X|}$ and $|\langle Z \setminus X \rangle| = \binom{n - |Z|}{k - |Z| + |X|}$.
\end{enumerate}
The analogous statements hold for $Y \in \cB_{\rm std}$ and the second mutation $\mu^2_Y(\cA, \cB)$.
\end{proposition}

\begin{proof}
\noindent\ref{itm:mutation-traces} Since $X \in \cA_{\rm std}$, every $k$-set containing $X$ belongs to $\cA$, hence $\langle X \rangle \subseteq \cA$. If $F \in \langle Z \setminus X \rangle \cap \cB$, then $F$ contains some $H \in \cB_{\rm std}$, and $H \subseteq F \cap Z = Z \setminus X$. Then $H \cap X = \emptyset$, which contradicts \cref{thm:generating}\ref{itm:cross-via-generating}.

\noindent\ref{itm:mutation-cross-intersecting} Since $\cA \setminus \langle X \rangle \subseteq \cA$, it suffices to show that $\langle Z \setminus X \rangle$ and $\cA \setminus \langle X \rangle$ are cross-intersecting. Let $F \in \langle Z \setminus X \rangle$ and $A \in \cA$ with $A \cap F = \emptyset$. Then $A$ contains some $G \in \cA_{\rm std}$, and $G \subseteq A \cap Z \subseteq Z \setminus (F \cap Z) = X$. Since $\cA_{\rm std}$ is an antichain, $G = X$, hence $X \subseteq A \cap Z \subseteq X$, that is, $A \in \langle X \rangle$.

\noindent\ref{itm:mutation-sizes} This follows from the formula for $|\langle X \rangle_S|$ in \cref{sec:operations}.

The last statement follows by symmetry.
\end{proof}

\subsection{Hinges} \label{sec:hinges}

To compare a pair with its mutations, we use two mutations in opposite directions, at members of the two standard generating families that together cover $Z$ and meet in only one element.

\begin{definition} \label{def:hinge}
A {\em hinge} of $(\cA, \cB)$ is a pair $(X, Y)$ with $X \in \cA_{\rm std}$, $Y \in \cB_{\rm std}$, $X \cup Y = Z$ and $|X \cap Y| = 1$. The mutations $\mu^1_X$ and $\mu^2_Y$ are called the {\em mutations at the hinge} $(X, Y)$.
\end{definition}

By \cref{thm:generating}\ref{itm:cross-via-generating}, every $X \in \cA_{\rm std}$ intersects every $Y \in \cB_{\rm std}$. Hence $|X| + |Y| \geq |X \cup Y| + 1$, and hinges are exactly the pairs of members of the two standard generating families that cover $Z$ with the smallest possible total size $|X| + |Y| = |Z| + 1$. Let $(X, Y)$ be a hinge with $X \cap Y = \{z\}$. Then $Z \setminus X = Y \setminus \{z\}$ and $Z \setminus Y = X \setminus \{z\}$, hence
\begin{equation} \label{eq:hinge-mutations}
\begin{aligned}
  \mu^1_X(\cA, \cB) & = \left(\cA \setminus \langle X \rangle, \ \cB \cup \langle Y \setminus \{z\} \rangle \right), \\
  \mu^2_Y(\cA, \cB) & = \left(\cA \cup \langle X \setminus \{z\} \rangle, \ \cB \setminus \langle Y \rangle \right).
\end{aligned}
\end{equation}

In general, a maximal cross-intersecting pair of nonempty families need not have a hinge. For example, if $k \geq 2$, the $k$-uniform families $\cA$ and $\cB$ generated by $\{\{1, 2\}, \{3, 4\}\}$ and $\{\{1, 3\}, \{1, 4\}, \{2, 3\}, \{2, 4\}\}$, respectively, form a maximal cross-intersecting pair, and $\cA_{\rm std}$ and $\cB_{\rm std}$ are these generating families. Hence $Z = [4]$, but $|X \cup Y| = 3$ for all $X \in \cA_{\rm std}$ and $Y \in \cB_{\rm std}$.

\subsection{Mutations at a hinge}

In \cref{sec:maximum}, we compare a pair with both mutations at a hinge, which requires the four families in the two mutations to be nonempty. \cref{prop:hinge-cases} shows that this holds unless the pair is distinguished.

\begin{proposition} \label{prop:hinge-cases}
Assume that $\cA$ and $\cB$ are nonempty, and let $(X, Y)$ be a hinge of $(\cA, \cB)$ with $X \cap Y = \{z\}$. Then the following statements hold.
\begin{enumerate}
\item \label{itm:hinge-y-singleton} If $|Y| = 1$, then $(\cA, \cB) = (\cA^Z, \cB^Z)$ and $|Z| \leq k$.
\item \label{itm:hinge-x-singleton} If $|X| = 1$, then $(\cA, \cB) = (\cB^Z, \cA^Z)$ and $|Z| \leq k$.
\item \label{itm:hinge-nonempty} If $|X| \geq 2$ and $|Y| \geq 2$, then the four families in $\mu^1_X(\cA, \cB)$ and $\mu^2_Y(\cA, \cB)$ are nonempty.
\end{enumerate}
\end{proposition}

\begin{proof}
\noindent\ref{itm:hinge-y-singleton} If $|Y| = 1$, then $Y = \{z\}$, and $X = Z$ since $X \cup Y = Z$. In particular, $|Z| = |X| \leq k$. Since every member of $\cA_{\rm std}$ is a subset of $Z$ and $\cA_{\rm std}$ is an antichain, $\cA_{\rm std} = \{Z\}$, hence $\cA = \cA^Z$. Since $n \geq 2 k$, a $k$-set intersects all $k$-sets containing $Z$ if and only if it intersects $Z$, hence $\cB = \cA^* = \cB^Z$.

\noindent\ref{itm:hinge-x-singleton} This follows from \ref{itm:hinge-y-singleton} by symmetry.

\noindent\ref{itm:hinge-nonempty} The families $\cA \cup \langle X \setminus \{z\} \rangle$ and $\cB \cup \langle Y \setminus \{z\} \rangle$ are clearly nonempty. Assume to the contrary that $\cA \setminus \langle X \rangle = \emptyset$. Then every member of $\cA$ contains $X$, hence every $k$-set containing $z$ intersects all members of $\cA$, that is, belongs to $\cA^* = \cB$. Thus $\{z\} \in \cB_{\rm std}^\uparrow$, so some member of the antichain $\cB_{\rm std}$ is contained in $\{z\} \subseteq Y$, which forces $Y \subseteq \{z\}$, contradicting $|Y| \geq 2$. Hence $\cA \setminus \langle X \rangle \neq \emptyset$, and by symmetry, $\cB \setminus \langle Y \rangle \neq \emptyset$.
\end{proof}

\section{Maximum pairs with a hinge} \label{sec:maximum}

In this section, we compare a maximal cross-intersecting pair with its two mutations at a hinge, and determine the maximum pairs with a hinge among the pairs whose first family is contained in a star. The size function enters the comparison only through \cref{lem:three-points}, so the same argument applies to every size function of finite cross-degree.

Applying mutations iteratively and enlarging the resulting families to maximal ones, one can increase the value of a size function until it is not smaller than its values at those mutations of the current pair that consist of nonempty families. The main result of this section, \cref{thm:hinge-comparison}, gives a strong necessary condition for this process to stop: the sizes of the two families are bounded in terms of the ranks of the mutations at a hinge. In particular, it gives a necessary condition for a cross-intersecting pair to be maximum, which is used in the proof of the main result.

Throughout this section, $U = [n]$ with $n \geq 2 k + 1$. For a maximal cross-intersecting pair $(\cA, \cB)$, we use the notation of \cref{sec:mutations}: $Z$ is the common support of $\cA_{\rm std}$ and $\cB_{\rm std}$, and $\langle X \rangle := \langle X \rangle_Z$ for $X \subseteq Z$.

For an integer $r \in [1, k]$, let
\[
  \theta_{n, k}(r) := \frac{k - r + 1}{n - 2 k} \binom{n - r - 1}{k - r + 1}.
\]
Note that $\theta_{n, k}(r + 1) = \frac{k - r}{n - r - 1} \theta_{n, k}(r) \lt \theta_{n, k}(r)$ for $r \in [1, k - 1]$. We now apply \cref{lem:three-points} to compare a maximal cross-intersecting pair with its mutations at a hinge $(X, Y)$. If $|X| = 1$ or $|Y| = 1$, then the pair is distinguished by \cref{prop:hinge-cases}. Otherwise, assume that the value of the size function does not increase under the two mutations at the hinge. The following theorem shows that $|\cA|$ and $|\cB|$ are then bounded in terms of the ranks of the first and the second mutation, respectively.

\begin{theorem} \label{thm:hinge-comparison}
Let $(\cA, \cB)$ be a maximal cross-intersecting pair of nonempty families in $\binom{[n]}{k}$, and let $(X, Y)$ be a hinge of $(\cA, \cB)$ with $|X| \geq 2$ and $|Y| \geq 2$. Let $f(x, y)$ be a function of cross-degrees at most $(c_1, c_2)$, where $c_i = 0$ or $c_i \geq 1$ for $i = 1, 2$. Assume that $f(|\cA|, |\cB|)$ is not smaller than the values of $f$ at the sizes of the two mutations at the hinge $(X, Y)$. Then
\begin{align*}
  |\cA| & \leq c_2 \ \theta_{n, k}(|X|) \leq c_2 \ \theta_{n, k}(2), \\
  |\cB| & \leq c_1 \ \theta_{n, k}(|Y|) \leq c_1 \ \theta_{n, k}(2).
\end{align*}
In particular, $|\cA| \lt \binom{n - 1}{k - 1}$ if $n \geq (k - 1)(c_2 + 2) + 3$, and $|\cB| \lt \binom{n - 1}{k - 1}$ if $n \geq (k - 1)(c_1 + 2) + 3$.
\end{theorem}

\begin{proof}
Let $u := |X|$, $v := |Y|$ and $m := n - |Z| = n - u - v + 1$. Let $z$ be the unique element of $X \cap Y$. Let
\[
  a_1 := |\langle X \setminus \{z\} \rangle|, \quad b_1 := |\langle Y \rangle|, \quad a_2 := |\langle X \rangle|, \quad b_2 := |\langle Y \setminus \{z\} \rangle|.
\]
By \cref{eq:hinge-mutations} and \cref{prop:mutation-properties}\ref{itm:mutation-traces}, the sizes of the first and second mutations at the hinge are $(|\cA| - a_2, |\cB| + b_2)$ and $(|\cA| + a_1, |\cB| - b_1)$, respectively. By \cref{prop:mutation-properties}\ref{itm:mutation-sizes},
\[
  a_1 = \binom{m}{k - u + 1}, \quad b_1 = \binom{m}{k - v}, \quad a_2 = \binom{m}{k - u}, \quad b_2 = \binom{m}{k - v + 1}.
\]
They are positive, since $u, v \leq k \leq n - k$. Moreover, $|\cA| \gt a_2$ and $|\cB| \gt b_1$, since the four families in the two mutations are nonempty by \cref{prop:hinge-cases}\ref{itm:hinge-nonempty}.

Since $n \geq 2 k + 1$, a direct calculation shows that
\[
  \frac{a_2}{a_1} = \frac{k - u + 1}{n - k - v + 1} \lt \frac{n - k - u + 1}{k - v + 1} = \frac{b_2}{b_1}.
\]
Therefore, \cref{lem:three-points} applies, and a direct calculation, using $n - k - v + 1 = m - (k - u)$, gives
\[
  |\cA| \leq c_2 \frac{n - k - v + 1}{n - 2 k} \binom{m}{k - u} = c_2 \frac{m}{n - 2 k} \binom{m - 1}{k - u},
\]
and
\[
  |\cB| \leq c_1 \frac{k - v + 1}{n - 2 k} \binom{m}{k - v + 1} = c_1 \frac{m}{n - 2 k} \binom{m - 1}{k - v}.
\]
Since $u \geq 2$ and $v \geq 2$, we have $m \leq n - u - 1$ and $m \leq n - v - 1$. Since $m \binom{m - 1}{j} = (j + 1) \binom{m}{j + 1}$ is nondecreasing in $m$, we obtain $|\cA| \leq c_2 \theta_{n, k}(u)$ and $|\cB| \leq c_1 \theta_{n, k}(v)$. The other inequalities follow, since $\theta_{n, k}$ is decreasing and $2 \leq u, v \leq k$.

If $n \geq (k - 1)(c + 2) + 3$ for some $c \geq 0$, then $n - 2 k \gt c (k - 1)$, hence
\[
  c \ \theta_{n, k}(2) = \frac{c (k - 1)}{n - 2 k} \binom{n - 3}{k - 1} \lt \binom{n - 3}{k - 1} \lt \binom{n - 1}{k - 1}.
\]
The last inequality holds since $k \geq |X| \geq 2$. Applying this to $c = c_1$ and $c = c_2$ gives the last statement.
\end{proof}

The proof of \cref{thm:main} reduces to pairs whose first family is contained in the star $\cA^{[1]} = \cB^{[1]}$, since a family contained in a star is intersecting and can replace several families at once, see \cref{lem:star-reduction}. The following corollary, in which only the constant $c_1$ matters, shows that such a pair with a hinge is distinguished if it is maximum among these pairs.

\begin{corollary} \label{cor:maximum-in-star}
Let $f(x, y)$ be a function of cross-degrees at most $(c_1, c_2)$, where $c_i = 0$ or $c_i \geq 1$ for $i = 1, 2$. Suppose that $n \geq (k - 1)(c_1 + 2) + 3$. Let $(\cA, \cB)$ be a maximal cross-intersecting pair of nonempty families in $\binom{[n]}{k}$ with $\cA \subseteq \cA^{[1]}$ and with a hinge. Assume that $f(|\cA|, |\cB|) \geq f(|\cA'|, |\cB'|)$ for every pair $(\cA', \cB')$ of nonempty cross-intersecting families in $\binom{[n]}{k}$ with $\cA' \subseteq \cA^{[1]}$. Then $(\cA, \cB) = (\cA^S, \cB^S)$ for some $S \subseteq [n]$ with $1 \in S$ and $|S| \leq k$.
\end{corollary}

\begin{proof}
Let $(X, Y)$ be a hinge of $(\cA, \cB)$ with $X \cap Y = \{z\}$. Every member $G$ of $\cA_{\rm std}$ contains $1$, since otherwise a $k$-set containing $G$ but not $1$ would belong to $\cA \subseteq \cA^{[1]}$. If $|Y| = 1$, then $(\cA, \cB) = (\cA^Z, \cB^Z)$ with $|Z| \leq k$ by \cref{prop:hinge-cases}\ref{itm:hinge-y-singleton}, and $1 \in Z$ since $Z = X \in \cA_{\rm std}$. If $|X| = 1$, then $\cB^Z = \cA \subseteq \cA^{[1]}$ by \cref{prop:hinge-cases}\ref{itm:hinge-x-singleton}, which forces $Z = \{1\}$, and $(\cA, \cB) = (\cA^{[1]}, \cB^{[1]})$.

Assume now that $|X| \geq 2$ and $|Y| \geq 2$. Since $\cB = \cA^* \supseteq (\cA^{[1]})^* = \cA^{[1]}$ (as $n \geq 2 k$), every $k$-set containing $1$ belongs to $\cB$, that is, $\{1\} \in \cB_{\rm std}^\uparrow$. Since $\cA$ is nonempty, $\cB \neq \binom{[n]}{k}$, hence $\emptyset \notin \cB_{\rm std}$ and $\{1\} \in \cB_{\rm std}$. Since $\cB_{\rm std}$ is an antichain and $|Y| \geq 2$, we have $1 \notin Y$. Together with $1 \in X$, this gives $1 \in X \setminus \{z\}$. Therefore, the first families $\cA \setminus \langle X \rangle$ and $\cA \cup \langle X \setminus \{z\} \rangle$ of the two mutations at the hinge are contained in $\cA^{[1]}$. By the assumption, together with \cref{prop:hinge-cases}\ref{itm:hinge-nonempty} and \cref{prop:mutation-properties}\ref{itm:mutation-cross-intersecting}, $f(|\cA|, |\cB|)$ is not smaller than the values of $f$ at the sizes of the two mutations. Thus, \cref{thm:hinge-comparison} implies $|\cB| \lt \binom{n - 1}{k - 1}$. On the other hand, $\cA^{[1]} \subseteq \cB$, hence $|\cB| \geq \binom{n - 1}{k - 1}$, which is a contradiction.
\end{proof}

\section{Distinguished cross-intersecting families} \label{sec:distinguished}

In this section, we study the distinguished pairs $(\cA^S, \cB^S)$. They form a chain of mutations between the pair of stars and $(\{[k]\}, \cB^{[k]})$, that is, between the two types of extremal structures discussed in the introduction. Comparing consecutive pairs in this chain gives the bounds on $|S|$ in \cref{thm:main}. Finally, we show that distinguished pairs are determined by their sizes. Throughout this section, $U = [n]$.

\subsection{The chain of distinguished pairs} \label{sec:chain}

Let $S \subseteq [n]$ with $s := |S| \in [1, k]$, and let $\cA^S$ and $\cB^S$ be as in \cref{thm:main}. For $S = \emptyset$, the same definitions give $(\cA^{\emptyset}, \cB^{\emptyset}) = (\binom{[n]}{k}, \emptyset)$. Note that $\cA^{\{i\}} = \cB^{\{i\}}$ is the star of all $k$-sets containing $i$. The sizes of these families are
\begin{equation} \label{eq:distinguished-sizes}
\begin{aligned}
  |\cA^S| & = \binom{n - s}{k - s}, \\
  |\cB^S| & = \binom{n}{k} - \binom{n - s}{k} = \sum_{j = 1}^s \binom{n - j}{k - 1}.
\end{aligned}
\end{equation}

Since $n \geq 2 k$, a $k$-set intersects all members of $\cA^S$ if and only if it intersects $S$, and a $k$-set intersects all members of $\cB^S$ if and only if it contains $S$. Hence $(\cA^S, \cB^S)$ is maximal cross-intersecting. The standard generating families are
\[
\cA^{S}_{\rm std} = \left\{ S \right \}
\quad \text{and} \quad
\cB^{S}_{\rm std} = \left\{ \{i\} \mid i \in S \right\},
\]
and the common support is $S$. The hinges of $(\cA^S, \cB^S)$ are exactly the pairs $(S, \{i\})$ with $i \in S$. The mutations at the hinge $(S, \{i\})$ are
\[
  \mu^1_{S}(\cA^S, \cB^S) = \left(\emptyset, \binom{[n]}{k}\right), \quad \text{and} \quad \mu^2_{\{i\}}(\cA^S, \cB^S) = (\cA^{S \setminus \{i\}}, \cB^{S \setminus \{i\}}).
\]
In particular, the second mutations of $(\cA^{[s]}, \cB^{[s]})$ form a chain
\[
  \left(\binom{[n]}{k}, \emptyset\right) = (\cA^{\emptyset}, \cB^{\emptyset}) \xleftarrow{\ \mu^2\ } (\cA^{[1]}, \cB^{[1]}) \xleftarrow{\ \mu^2\ } \cdots \xleftarrow{\ \mu^2\ } (\cA^{[k]}, \cB^{[k]}),
\]
where each arrow is the second mutation at $\{s\} \in \cB^{[s]}_{\rm std}$.

\subsection{Comparing distinguished families} \label{sec:comparing}

We now show that $(\cA^{[s]}, \cB^{[s]})$ is not maximum unless $s$ is small or $s = k$. To this end, we compare three consecutive pairs in the chain of mutations in \cref{sec:chain}:
\[
    (\cA^{[s - 1]}, \cB^{[s - 1]}) \xleftarrow{\ \mu^2\ } (\cA^{[s]}, \cB^{[s]}) \xleftarrow{\ \mu^2\ } (\cA^{[s + 1]}, \cB^{[s + 1]}).
\]

\begin{proposition} \label{prop:consecutive}
Let $f(x, y)$ be a function of cross-degrees at most $(c_1, c_2)$, where $c_i = 0$ or $c_i \geq 1$ for $i = 1, 2$. Suppose that $n \geq 2 k + 1$. Let $s \in [2, k - 1]$ be an integer such that
\begin{align*}
  f(|\cA^{[s]}|, |\cB^{[s]}|) & \geq f(|\cA^{[s - 1]}|, |\cB^{[s - 1]}|), \quad \text{and} \\
  f(|\cA^{[s]}|, |\cB^{[s]}|) & \geq f(|\cA^{[s + 1]}|, |\cB^{[s + 1]}|).
\end{align*}
Then
\begin{equation} \label{eq:consecutive}
  \left(\frac{s}{s - 1}\right)^{c_1} \gt 1 + \frac{(n - 2 k)(n - s)}{(k - s + 1)(2 n - 2 s - k + 1)}.
\end{equation}
In particular,
\[
  s \lt \frac{c_1 n}{n - 2 k + c_1} \quad \text{and} \quad n \lt (k - 1) 2^{c_1 + 1} + 2.
\]
If moreover $n \geq (c_1 + 1)(2 k - c_1)$, then $s \leq \lceil c_1 \rceil$.
\end{proposition}

\begin{proof}
Let
\begin{align*}
a_1 & := |\cA^{[s - 1]}| - |\cA^{[s]}| = \binom{n - s}{k - s + 1}, & b_1 & := |\cB^{[s]}| - |\cB^{[s - 1]}| = \binom{n - s}{k - 1}, \\
a_2 & := |\cA^{[s]}| - |\cA^{[s + 1]}| = \binom{n - s - 1}{k - s}, & b_2 & := |\cB^{[s + 1]}| - |\cB^{[s]}| = \binom{n - s - 1}{k - 1},
\end{align*}
which are positive, and let $a := |\cA^{[s]}|$ and $b := |\cB^{[s]}|$. Then the assumptions say that $f(a, b) \geq f(a + a_1, b - b_1)$ and $f(a, b) \geq f(a - a_2, b + b_2)$. A direct calculation shows that
\[
  \frac{a_2}{a_1} = \frac{k - s + 1}{n - s} \lt \frac{n - s - k + 1}{n - s} = \frac{b_2}{b_1},
\]
since $n \geq 2 k + 1$. Moreover, $a \gt a_2$ and $b \gt b_1$, since $a - a_2 = |\cA^{[s + 1]}|$ and $b - b_1 = |\cB^{[s - 1]}|$ are positive, so \cref{lem:three-points} applies. In particular, $c_1 \geq 1$, and \cref{eq:ratio-b} gives
\[
  \left(\frac{b}{b - b_1}\right)^{c_1} \geq 1 + q,
\]
where
\[
  q := \frac{a_1 b_2 - a_2 b_1}{a_2 (b_1 + b_2)} = \frac{(n - 2 k)(n - s)}{(k - s + 1)(2 n - 2 s - k + 1)}.
\]
Since
\[
  b - b_1 = |\cB^{[s - 1]}| = \sum_{j = 1}^{s - 1} \binom{n - j}{k - 1} \gt (s - 1) \binom{n - s}{k - 1} = (s - 1) b_1,
\]
we have $\frac{b}{b - b_1} = 1 + \frac{b_1}{b - b_1} \lt \frac{s}{s - 1}$, which proves \cref{eq:consecutive}.

By Bernoulli's inequality and \cref{eq:consecutive}, $1 - \frac{c_1}{s} \leq (1 - \frac{1}{s})^{c_1} \lt \frac{1}{1 + q}$, hence
\[
  s \lt c_1 \frac{1 + q}{q} = c_1 \left(1 + \frac{(k - s + 1)(2 n - 2 s - k + 1)}{(n - 2 k)(n - s)}\right) \leq c_1 \frac{n - s}{n - 2 k},
\]
where the last inequality holds since $s \geq 2$. Hence $s (n - 2 k + c_1) \lt c_1 n$. If $n \geq (c_1 + 1)(2 k - c_1)$, then $c_1 n \leq (c_1 + 1)(n - 2 k + c_1)$, hence $s \lt c_1 + 1$, that is, $s \leq \lceil c_1 \rceil$.

Finally, since $k - s + 1 \leq k - 1$, $2 n - 2 s - k + 1 \leq 2 (n - s)$ and $\frac{s}{s - 1} \leq 2$, \cref{eq:consecutive} gives $2^{c_1} \gt 1 + q \geq 1 + \frac{n - 2 k}{2 (k - 1)} = \frac{n - 2}{2 (k - 1)}$, that is, $n \lt (k - 1) 2^{c_1 + 1} + 2$.
\end{proof}

\subsection{Distinguished families are determined by their sizes} \label{sec:determined-by-sizes}

In this subsection, we show that a cross-intersecting pair with the same sizes as $(\cA^{[s]}, \cB^{[s]})$ is distinguished, see \cref{prop:determined-by-sizes}. This is used in \cref{sec:proof} to recover maximum families from their sizes. \cref{prop:determined-by-sizes} is a reformulation of a special case of \cite[Proposition~2.3]{SFQ22}. Its equality case also follows from results of F\"uredi--Griggs \cite{FG86} and M\"ors \cite{Mor85}. We give a short proof using the following version of the Kruskal--Katona theorem due to Lov\'asz \cite{Lov79}, together with its equality case.

For an integer $j \geq 1$, the function $t \mapsto \binom{t}{j}$ of a real variable $t$ is strictly increasing on $[j - 1, \infty)$ and vanishes at $j - 1$. The {\em shadow} of a family $\cF$ of $r$-sets is
\[
  \partial \cF := \left\{ G \mid |G| = r - 1, \ G \subseteq F \text{ for some } F \in \cF \right\}.
\]

\begin{theorem} [Lov\'asz {\cite[Problem~13.31(b)]{Lov79}}] \label{thm:lovasz}
Let $r \geq 1$, and let $\cF \subseteq \binom{[n]}{r}$ be a family with $|\cF| = \binom{t}{r}$ for some real number $t \geq r$. Then $|\partial \cF| \geq \binom{t}{r - 1}$. If moreover $t$ is an integer and $|\partial \cF| = \binom{t}{r - 1}$, then $\cF = \binom{M}{r}$ for some $M \subseteq [n]$ with $|M| = t$.
\end{theorem}

The following proposition is proved by applying \cref{thm:lovasz} to the complements of the members of one of the two families and to their iterated shadows.

\begin{proposition} \label{prop:determined-by-sizes}
Assume that $n \geq 2 k + 1$. Let $\cA, \cB \subseteq \binom{[n]}{k}$ be two cross-intersecting families with $|\cA| = |\cA^{[s]}|$ for some $s \in [k]$. Then $|\cB| \leq |\cB^{[s]}|$. Moreover, equality holds if and only if $(\cA, \cB) = (\cA^T, \cB^T)$ for some $T \subseteq [n]$ with $|T| = s$.
\end{proposition}

\begin{proof}
Let $\cC_0 := \{ [n] \setminus A \mid A \in \cA \}$, and $\cC_i := \partial \cC_{i - 1}$ for $i \in [n - 2 k]$. Then $\cC_{n - 2 k}$ is the family of $k$-sets that are disjoint from some member of $\cA$, hence it is disjoint from $\cB$. For $i \in [0, n - 2 k]$, let $t_i \geq n - k - i$ be the real number with $|\cC_i| = \binom{t_i}{n - k - i}$, which exists since $\cC_i \neq \emptyset$. Since $|\cC_0| = \binom{n - s}{k - s} = \binom{n - s}{n - k}$, we have $t_0 = n - s$, and \cref{thm:lovasz} gives $t_0 \leq t_1 \leq \dots \leq t_{n - 2 k}$. Hence $|\cC_{n - 2 k}| \geq \binom{n - s}{k}$, and $|\cB| \leq \binom{n}{k} - \binom{n - s}{k} = |\cB^{[s]}|$ by \cref{eq:distinguished-sizes}.

If $|\cB| = |\cB^{[s]}|$, then $t_{n - 2 k} \leq n - s$, hence $t_1 = t_0 = n - s$. By \cref{thm:lovasz}, $\cC_0 = \binom{M}{n - k}$ for some $M \subseteq [n]$ with $|M| = n - s$. Let $T := [n] \setminus M$, so that $\cA = \cA^T$. Since $\cA$ and $\cB$ are cross-intersecting, $\cB \subseteq (\cA^T)^* = \cB^T$, see \cref{sec:chain}, and $\cB = \cB^T$ since $|\cB| = |\cB^T|$. Conversely, $|\cB^T| = |\cB^{[s]}|$ by \cref{eq:distinguished-sizes}.
\end{proof}

\section{Hinges of left-shifted families} \label{sec:left-shifted}

In this section, we show that every maximal cross-intersecting pair of nonempty left-shifted families has a hinge, see \cref{prop:hinge-exists}. This is the only step of our proof that uses shifting. In \cref{sec:proof}, left-shifted families arise as initial segments of the lexicographic order. Throughout this section, $U = [n]$.

For $i, j \in U$, define the {\em shifting operator} $s_{ij} \from 2^{U} \to 2^{U}$ by
\[
F \mapsto \begin{cases}
(F \setminus \{j\}) \cup \{i\}, & i \notin F, j \in F, \\
F, & \text{otherwise}.
\end{cases}
\]
A family $\cF \subseteq 2^U$ is {\em $s_{ij}$-stable} if $s_{ij}(F) \in \cF$ for every $F \in \cF$. It is {\em left-shifted} if it is $s_{ij}$-stable for all $1 \leq i \lt j \leq n$, see for example \cite{Fra87}.

To find a hinge, we use two properties of the standard generating family of a left-shifted family: its up-closure is again left-shifted, and its support is of the form $[\ell]$.

\begin{lemma} \label{lem:std-shifts}
Let $\cF \subseteq \binom{U}{k}$ be a family, and $\cF_{\rm std}$ its standard generating family. Then the following statements hold.
\begin{enumerate}
  \item \label{itm:std-stable} If $\cF$ is $s_{ij}$-stable for some $i, j \in U$, then $\cF_{\rm std}^\uparrow$ is also $s_{ij}$-stable.
  \item \label{itm:std-support-interval} If $\cF$ is left-shifted, then $\supp \cF_{\rm std} = [\ell]$ for some $\ell \in [0, n]$.
\end{enumerate}
\end{lemma}

\begin{proof}
\noindent\ref{itm:std-stable} Since $\cF$ is $s_{ij}$-stable, $\cF^\uparrow$ is also $s_{ij}$-stable. Indeed, if $W \in \cF^\uparrow$ with $i \notin W$ and $j \in W$, take $F \in \cF$ with $F \subseteq W$. If $j \in F$, then $s_{ij}(F) \in \cF$ and $s_{ij}(F) \subseteq s_{ij}(W)$, and otherwise $F \subseteq s_{ij}(W)$. Either way, $s_{ij}(W) \in \cF^\uparrow$.

Now let $G \in \cF_{\rm std}^\uparrow$ with $s_{ij}(G) \neq G$, since otherwise there is nothing to prove. Then $i \notin G$ and $j \in G$. Let $X \in s_{ij}(G)^\uparrow \cap \binom{U}{\geq k}$ be arbitrary. If $j \in X$, then $G \subseteq X$, hence $X \in \cF^\uparrow$. If $j \notin X$, then $G \subseteq s_{ji}(X)$ and $|s_{ji}(X)| = |X| \geq k$, hence $s_{ji}(X) \in \cF^\uparrow$ since $G \in \cF_{\rm std}^\uparrow$. Therefore, $X = s_{ij}(s_{ji}(X)) \in \cF^\uparrow$ since $\cF^\uparrow$ is $s_{ij}$-stable. Since $X$ is arbitrary, $s_{ij}(G) \in \cF_{\rm std}^\uparrow$.

\noindent\ref{itm:std-support-interval} Let $j \in \supp \cF_{\rm std}$ and $1 \leq i \lt j$. We need to show that $i \in \supp \cF_{\rm std}$. Take $G \in \cF_{\rm std}$ with $j \in G$, and assume that $i \notin G$, since otherwise we are done. Then $s_{ij}(G) \in \cF_{\rm std}^\uparrow$ by \ref{itm:std-stable}, hence there exists $G' \in \cF_{\rm std}$ with $G' \subseteq s_{ij}(G) = (G \setminus \{j\}) \cup \{i\}$. If $i \notin G'$, then $G' \subsetneq G$, contradicting the fact that $\cF_{\rm std}$ is an antichain. Hence $i \in G' \subseteq \supp \cF_{\rm std}$.
\end{proof}

A maximal cross-intersecting pair of nonempty families need not have a hinge, see \cref{sec:hinges}. The following proposition is where left-shifted families are needed: they guarantee that hinges exist. For a maximal cross-intersecting pair $(\cA, \cB)$, we use the notation of \cref{sec:mutations}: $Z$ is the common support of $\cA_{\rm std}$ and $\cB_{\rm std}$, and $\langle X \rangle := \langle X \rangle_Z$ for $X \subseteq Z$.

\begin{proposition} \label{prop:hinge-exists}
Let $(\cA, \cB)$ be a maximal cross-intersecting pair of nonempty left-shifted families in $\binom{[n]}{k}$. Then $(\cA, \cB)$ has a hinge.
\end{proposition}

\begin{proof}
By \cref{lem:std-shifts}\ref{itm:std-support-interval} and \cref{thm:generating}\ref{itm:common-support}, $Z = [\ell]$ for some $\ell \in [0, n]$. If $\ell = 0$, then $\cA_{\rm std} \subseteq \{\emptyset\}$, and $\cA_{\rm std} = \{\emptyset\}$ since $\cA$ is nonempty. Then $\cA = \binom{[n]}{k}$, and $\cB = \cA^* = \emptyset$ since $n \geq 2 k$, which is a contradiction. Hence $\ell \geq 1$.

Let $X \in \cA_{\rm std}$ with $\ell \in X$, which exists since $\ell \in \supp \cA_{\rm std}$, and let $Y := ([\ell] \setminus X) \cup \{\ell\}$. We show that $(X, Y)$ is a hinge. First, we find $H \in \cB_{\rm std}$ with $\ell \in H \subseteq Y$. By the minimality of $X$, we have $X \setminus \{\ell\} \notin \cA_{\rm std}^\uparrow$, hence some $k$-set $F$ containing $X \setminus \{\ell\}$ does not belong to $\cA$. Since $\cA = \cB^*$, there exists $B \in \cB$ with $F \cap B = \emptyset$. Take $H \in \cB_{\rm std}$ with $H \subseteq B$. Then $H \subseteq B \cap Z \subseteq Z \setminus (X \setminus \{\ell\}) = Y$. By \cref{thm:generating}\ref{itm:cross-via-generating}, $H \cap X \neq \emptyset$, hence $\ell \in H$.

We claim that $H = Y$. Otherwise, take $w \in Y \setminus H$, and note that $w \lt \ell$ and $w \notin X$. By \cref{lem:std-shifts}\ref{itm:std-stable}, $\cB_{\rm std}^\uparrow$ is $s_{w \ell}$-stable, hence $(H \setminus \{\ell\}) \cup \{w\} \in \cB_{\rm std}^\uparrow$. However, $(H \setminus \{\ell\}) \cup \{w\} \subseteq Y \setminus \{\ell\}$ is disjoint from $X$, so $\cB_{\rm std}$ contains a member disjoint from $X \in \cA_{\rm std}$, contradicting \cref{thm:generating}\ref{itm:cross-via-generating}. Hence $Y = H \in \cB_{\rm std}$, and $(X, Y)$ is a hinge.
\end{proof}

\section{Proof of the main theorem} \label{sec:proof}

\cref{cor:maximum-in-star} determines the maximum pairs with a hinge among the pairs whose first family is contained in a star. We now use Hilton's lemma to obtain hinges, and reduce \cref{thm:main} to such pairs. Note that the assumption $n \geq (k - 1)(c + 2) + 3$ in \cref{thm:main} implies $n \geq 2 k + 1$.

We use the {\em lexicographic order} on $\binom{[n]}{k}$: for distinct $F, F' \in \binom{[n]}{k}$, we write $F \lt_{\rm lex} F'$ if $\min (F \setminus F') \lt \min (F' \setminus F)$. Let $\cL(n, k, m)$ denote the initial segment of length $m$ with respect to the lexicographic order on $\binom{[n]}{k}$, that is, the family of the first $m$ sets. Initial segments are totally ordered by inclusion, and they are left-shifted: if $F \in \cL(n, k, m)$ and $i \lt j$ with $i \notin F$ and $j \in F$, then $s_{ij}(F) \lt_{\rm lex} F$, hence $s_{ij}(F) \in \cL(n, k, m)$. Moreover, for $s \in [k]$, the $k$-sets containing $[s]$ precede all other $k$-sets, hence
\[
  \cL\left(n, k, \binom{n - s}{k - s}\right) = \cA^{[s]}.
\]
In particular, $\cL(n, k, \binom{n - 1}{k - 1}) = \cA^{[1]}$ is a star.

In the study of cross-intersecting families, a key tool is Hilton's lemma, which follows from the Kruskal--Katona theorem \cite{Kat68, Kru63}, see also \cite[Lemma 152]{GP18}.

\begin{lemma} [Hilton \cite{Hil76}] \label{lem:hilton}
Let $\cF, \cF'$ be two families in $\binom{[n]}{k}$. If $\cF$ and $\cF'$ are cross-intersecting, then $\cL(n, k, |\cF|)$ and $\cL(n, k, |\cF'|)$ are also cross-intersecting.
\end{lemma}

Hilton's lemma also implies that two cross-intersecting families cannot both have more members than a star.

\begin{corollary} \label{cor:one-large}
Assume that $n \geq 2 k$. Let $\cF, \cF' \subseteq \binom{[n]}{k}$ be cross-intersecting families. Then at most one of $|\cF|$ and $|\cF'|$ is larger than $\binom{n - 1}{k - 1}$.
\end{corollary}

\begin{proof}
An initial segment of length larger than $\binom{n - 1}{k - 1}$ contains the star $\cA^{[1]}$ and a $k$-set $F$ with $1 \notin F$. Since $n \geq 2 k$, some member of $\cA^{[1]}$ is disjoint from $F$. Hence two such initial segments are not cross-intersecting, and the claim follows from \cref{lem:hilton}.
\end{proof}

For the rest of this section, we fix a function $f$ of $d$ variables with cross-degree at most $c$, where $c = 0$ or $c \geq 1$. For $p \in [d]$, let
\[
  g_p(x, y) := f(x, \dots, x, y, x, \dots, x),
\]
where $y$ is in the $p$-th position. By \cref{lem:two-variables}, $g_p$ has cross-degrees at most $(c, (d - 1) c)$. Since $c = 0$ or $c \geq 1$, the same holds for $(d - 1) c$. In the following, only the constant $c_1 = c$ of $g_p$ matters, so the bounds in \cref{thm:main} do not depend on $d$.

We say that nonempty cross-intersecting families $\cF_1, \dots, \cF_d \subseteq \binom{[n]}{k}$ are {\em maximum} if $f(|\cF_1|, \dots, |\cF_d|)$ is maximum among all nonempty cross-intersecting families in $\binom{[n]}{k}$. Since $f$ is strictly increasing in all variables, maximum families are maximal cross-intersecting. In particular, if $(\cF_1, \dots, \cF_d)$ is maximum, then $\cF_i = (\bigcup_{j \neq i} \cF_j)^*$ for every $i$.

The following lemma treats the case in which $d - 1$ of the families are contained in the star $\cA^{[1]}$. These $d - 1$ families are then equal, and \cref{cor:maximum-in-star} applies to the pair formed by one of them and the remaining family.

\begin{lemma} \label{lem:star-reduction}
Suppose that $n \geq (k - 1)(c + 2) + 3$. Let $(\cG_1, \dots, \cG_d)$ be maximum nonempty left-shifted cross-intersecting families in $\binom{[n]}{k}$, at most one of which is not contained in $\cA^{[1]}$. Then there exists a nonempty $S \subseteq [n]$ with $|S| \leq k$ such that $(\cG_1, \dots, \cG_d)$ is a permutation of $(\cA^S, \dots, \cA^S, \cB^S)$.
\end{lemma}

\begin{proof}
Take $p \in [d]$ such that $\cG_j \subseteq \cA^{[1]}$ for all $j \neq p$, and let $\cA := \bigcup_{j \neq p} \cG_j$ and $\cB := \cG_p$. Then $\cA \subseteq \cA^{[1]}$ is intersecting, and $\cA$ and $\cB$ are cross-intersecting. Therefore, replacing $\cG_j$ by $\cA$ for all $j \neq p$, we obtain cross-intersecting families. By maximality, $\cG_j = \cA$ for all $j \neq p$, and $\cB = \cA^*$. Since $\cA \subseteq \cA^{[1]}$, we have $\cB = \cA^* \supseteq (\cA^{[1]})^* = \cA^{[1]}$ (as $n \geq 2 k$), hence $\cB^* \subseteq \cA^{[1]}$ is intersecting. Replacing $\cG_j$ by $\cB^* \supseteq \cA$ for all $j \neq p$, we again obtain cross-intersecting families, hence $\cA = \cB^*$ by maximality. Therefore, $(\cA, \cB)$ is a maximal cross-intersecting pair of nonempty left-shifted families with $\cA \subseteq \cA^{[1]}$, and it has a hinge by \cref{prop:hinge-exists}.

For every pair $(\cA', \cB')$ of nonempty cross-intersecting families with $\cA' \subseteq \cA^{[1]}$, the family $\cA'$ is intersecting. Hence the $d$ families with $\cA'$ at the positions in $[d] \setminus \{p\}$ and $\cB'$ at the position $p$ are cross-intersecting, and $g_p(|\cA'|, |\cB'|) \leq f(|\cG_1|, \dots, |\cG_d|) = g_p(|\cA|, |\cB|)$. By \cref{cor:maximum-in-star} applied to $g_p$ with $c_1 = c$, we have $(\cA, \cB) = (\cA^S, \cB^S)$ for some nonempty $S \subseteq [n]$ with $|S| \leq k$.
\end{proof}

To apply \cref{lem:star-reduction}, we replace the families by initial segments of the lexicographic order of the same sizes, all but one of which are contained in the star $\cA^{[1]}$ by \cref{cor:one-large}.

\begin{proof}[Proof of \cref{thm:main}]
Let $\cF_1, \dots, \cF_d$ be maximum nonempty cross-intersecting families. For $j \in [d]$, let $\cL_j := \cL(n, k, |\cF_j|)$. Then the families $\cL_1, \dots, \cL_d$ are nonempty and left-shifted, they are cross-intersecting by \cref{lem:hilton}, and they are also maximum. By \cref{cor:one-large}, at most one of them has more than $\binom{n - 1}{k - 1}$ members, and the others are contained in $\cL(n, k, \binom{n - 1}{k - 1}) = \cA^{[1]}$. By \cref{lem:star-reduction}, there exist a nonempty $S' \subseteq [n]$ with $s := |S'| \leq k$ and $p \in [d]$ such that $\cL_p = \cB^{S'}$ and $\cL_j = \cA^{S'}$ for $j \neq p$. Therefore $|\cF_p| = |\cB^{[s]}|$ and $|\cF_j| = |\cA^{[s]}|$ for $j \neq p$.

For $j \neq p$, \cref{prop:determined-by-sizes} applied to $(\cF_j, \cF_p)$ gives $T_j \subseteq [n]$ with $|T_j| = s$ and $(\cF_j, \cF_p) = (\cA^{T_j}, \cB^{T_j})$. Since $n - s \geq k$, the complement $[n] \setminus T_j$ is the union of all $k$-subsets of $[n]$ that are not in $\cB^{T_j} = \cF_p$. Hence all the sets $T_j$ with $j \neq p$ are the same set $S$, and $(\cF_1, \dots, \cF_d)$ is a permutation of $(\cA^S, \dots, \cA^S, \cB^S)$.

Finally, assume that $s \in [2, k - 1]$. Since the families with $\cA^{[s \pm 1]}$ at the positions in $[d] \setminus \{p\}$ and $\cB^{[s \pm 1]}$ at the position $p$ are nonempty and cross-intersecting, the value $g_p(|\cA^{[s]}|, |\cB^{[s]}|)$ is not smaller than $g_p(|\cA^{[s \pm 1]}|, |\cB^{[s \pm 1]}|)$. By \cref{prop:consecutive} applied to $g_p$ with $c_1 = c$, we have $n \lt (k - 1) 2^{c + 1} + 2$, which proves \cref{thm:main}\ref{itm:main-one-or-k}. Moreover, $s \leq \lceil c \rceil$ if $n \geq (c + 1)(2 k - c)$, in particular if $n \geq 2 (c + 1) k$. This proves \cref{thm:main}\ref{itm:main-at-most-c}.
\end{proof}

\section*{Acknowledgments}

The first author is supported by NSFC (12350710787) and NSFC (12471311), and the second author is supported by NSFC (12371332).

\section*{AI disclosure}

The initial unification (the setup in Example~\ref{ex:interpolation}) was obtained in March 2023, and the main result was obtained in July 2025; neither involved any use of AI. AI tools were used only from October 2026 through the final writing stage.

\end{document}